\documentclass[12pt]{article}
\usepackage{amsfonts,amsmath,amssymb,amscd}
\usepackage{shadow}
\usepackage{graphicx}
\usepackage{graphicx}
\usepackage{color}

\newtheorem{theo}{Theorem}[section]

\makeatletter 
\@addtoreset{equation}{section}
\makeatother   

\def\qed{\hfill \rule{4pt}{7pt}}

\def\pf{\noindent {\it Proof.} }

\begin{document}

\begin{center}
{\Large \bf  Disposition Polynomials and Plane Trees}
\end{center}

\begin{center}
William Y. C. Chen$^{1}$and Janet F.F. Peng$^{2}$\\[6pt]

$^{1,2}$Center for Combinatorics, LPMC-TJKLC\\
Nankai University, Tianjin 300071, P. R. China

 Email: $^{1}${\tt
chen@nankai.edu.cn}, $^{2}${\tt janet@mail.nankai.edu.cn}
\end{center}

\vskip 6mm \noindent {\bf Abstract.} We define the disposition
polynomial $R_{m}(x_1, x_2, \ldots, x_n)$ as
$\prod_{k=0}^{m-1}(x_1+x_2+\cdots+x_n+k)$. When $m=n-1$, this
polynomial becomes  the generating function of plane trees with
respect to certain statistics as given by Guo and Zeng. When $x_i=1$
for $1\leq i\leq n$,
 $R_{m}(x_1, x_2, \ldots, x_n)$
reduces to the rising factorial $n(n+1)\cdots (n+m-1)$. Guo and Zeng
asked the question of finding a combinatorial proof of the formula
for the generating function of plane trees with respect to the
number of younger children and the number of elder children. We find
a combinatorial interpretation of the disposition polynomials in
terms of the number of right-to-left minima of each linear order in
a disposition. Then we establish a bijection between plane
trees on $n$ vertices
 and dispositions from $\{1, 2,\ldots, n-1\}$
to $\{ 1, 2,\ldots, n\}$ in the spirit of the Pr\"ufer correspondence.
  It gives an answer to the
question of Guo and Zeng, and it also
provides an answer to another question of Guo and Zeng
concerning an identity on the plane tree
expansion of a polynomial introduced by Gessel and Seo.



\section{Introduction}

The notation of dispositions was   introduced by Mullin and Rota
\cite{Rota1}, see also, Joni, Rota and Sagan \cite{Rota}.
Let $(x)^n$ denote the rising factorial $x(x+1)\cdots(x+n-1)$.
 Assume that $x$ is a nonnegative integer.
 Then $(x)^n$ can be interpreted
 as the number of dispositions from $[n]=\{1, 2, \ldots, n\}$ to $[x]=
 \{1, 2, \ldots, x\}$, where a disposition from $[n]$ to
 $[x]$ is a function from $[n]$ to $[x]$
 in which the pre-images  of each $i\in [x]$
 are endowed with a linear order. In other words,
 a disposition from $[n]$ to $[x]$ can be viewed as
 a decomposition of a permutation of $[n]$ into $x$
 parts.

In this paper, we introduce the disposition polynomials $R_{m}(x_1,
x_2, \ldots, x_n)$ as a multivariate extension of the rising
factorials by considering the number of right-to-left minima of
each linear order in a disposition from $[m]$ to $[n]$. More
 precisely, the disposition polynomials are defined by
\begin{equation}
R_{m}(x_1, x_2, \ldots,
x_n)=\prod_{k=0}^{m-1}(x_1+x_2+\cdots+x_n+k).
\end{equation}

As will be seen, for the disposition polynomials $R_{m}(x_1, x_2,
\ldots, x_n)$, the exponent of $x_i$ records the number of
right-to-left minima of the $i$-th linear order in a disposition.
For the purpose of this paper, we shall use the homogeneous disposition polynomials as given by
\begin{equation}
Q_{m}(x_1, x_2, \ldots,
x_n,t)=\prod_{k=0}^{m-1}(x_1+x_2+\cdots+x_n+kt).
\end{equation}
Note that the homogenous disposition polynomials have essentially the
same combinatorial interpretation as that of the disposition polynomials.

In fact, we are led to the above definition of the disposition polynomials
by the special case $m=n-1$ given by  Guo and Zeng \cite{GZ} for the enumeration of plane trees.

Let $\mathcal{P}_n$ denote the set of plane trees on $[n]$, where a plane tree on $[n]$ is a labeled rooted tree on $[n]$
 in which the children of each
vertex are linearly ordered,  and let $\mathcal{P}_n^{(r)}$ denote the set of plane trees on $[n]$ with root $r$. For $T\in\mathcal{P}_n$, let $i$ be a vertex of $T$ and $j$
be a child of $i$. If the smallest descendent of $j$ is smaller than  those of its brothers on the right, then $j$ is called a younger child of $i$. Otherwise, $j$ is called  an elder child of $i$,
that is, the smallest descendent of $j$ is bigger than those of a brother on the right. Denote by ${\rm young}_{T}(i)$ the number of younger children of $i$ in $T$, and denote by ${\rm eld}(T)$ the total
 number of elder vertices in $T$.  Guo and Zeng \cite{GZ} obtained the
 following formulas
\begin{equation}\label{e10}
\sum_{T\in \mathcal{P}_n}t^{{\rm eld}(T)}\prod_{i=1}^{n}x_i^{{\rm
young}_T(i)}=\prod_{k=0}^{n-2}(x_1+x_2+\cdots+x_n+kt),
\end{equation}
and
\begin{equation}\label{e11}
\sum_{T\in \mathcal{P}_n^{(r)}}t^{{\rm eld}(T)}\prod_{i=1}^{n}x_i^{{\rm
young}_T(i)}=x_r\prod_{k=1}^{n-2}(x_1+x_2+\cdots+x_n+kt).
\end{equation}

Guo and Zeng  proved the above formulas (\ref{e10}) and (\ref{e11}) by
induction and asked for combinatorial proofs. In answer to the questions of
Guo and Zeng, we first give a combinatorial interpretation of the
disposition polynomials. Then, for the case $m=n-1$,
we establish a Pr\"ufer type
 correspondence between plane trees and dispositions, which
implies combinatorial interpretations of both relations (\ref{e10})
and (\ref{e11}).

Replacing $n$ by $n+1$, $t$ by $t-z$ and
setting $r=1$, $x_1=x$ and $x_i=z$ for $2\leq i\leq n+1$,
the right hand side of (\ref{e11}) becomes the polynomial
\[x\prod_{k=1}^{n-1}(x+(n-k)z+kt),\]
which is the polynomial $P_n(t,z,x)$  introduced by Gessel and Seo \cite{Gessel}
 for the enumeration of labeled trees by the number of proper vertices.
Several
 expansions of the polynomial
$P_n(t,z,x)$ have been given by Gessel and Seo
\cite{Gessel}
in terms of rooted trees with proper vertices, $k$-ary trees with
proper vertices, $k$-colored ordered forests with proper vertices
and parking functions with lucky cars by using generating functions. Combinatorial proofs of
some of these relations have been found by Seo and Shin \cite{Shin} and by Shin
\cite{Shin1}.

With the above substitutions, (\ref{e11}) reduces to the
relation
\begin{equation}\label{e14}
\sum_{T \in
\mathcal{P}_{n+1}^{(1)}}x^{{\rm{young}}_{T}(1)}(t-z)^{{\rm{eld}}(T)}
z^{n-{\rm{young}}_{T}(1)-{\rm{eld}}(T)}=x\prod_{k=1}^{n-1}(x+(n-k)z+kt).
\end{equation}
Guo and Zeng \cite{GZ}  deduced the above formula as another combinatorial interpretation of the polynomial $P_n(t,z,x)$ of Gessel and Seo, and they raised the
question of finding a combinatorial interpretation of (\ref{e14}).

Our correspondence between plane trees and dispositions
can be directly applied to give a combinatorial interpretation of
(\ref{e14}).
Indeed,  the above relation
holds for plane trees with any given root $r$, that is,
\begin{equation}\label{ei15}
\sum_{T\in
\mathcal{P}_{n+1}^{(r)}}x^{{\rm{young}}_{T}(r)}(t-z)^{{\rm{eld}}(T)}
z^{n-{\rm{young}}_{T}(r)-{\rm{eld}}(T)}=x\prod_{k=1}^{n-1}(x+(n-k)z+kt).
\end{equation}

This paper is organized as follows. In Section 2, we  give a combinatorial explanation of the disposition polynomials.
Section 3 provides a Pr\"ufer type correspondence between plane trees and
dispositions which leads to combinatorial interpretations of ($\ref{e10}$) and ($\ref{e11}$).
Section $4$ is devoted to the combinatorial proof of ($\ref{ei15}$).

\section{The generating function of dispositions}

In this section, we give a combinatorial interpretation
of the disposition polynomials
\[ R_{m}(x_1, x_2, \ldots, x_n)=   \prod_{k=0}^{m-1}(x_1+x_2+\cdots+x_n+k). \]
The notion of dispositions
was introduced by Mullin and Rota \cite{Rota1} as a combinatorial
explanation of the rising factorials $(x)^n=x(x+1)\cdots (x+n-1)$,
see also Joni, Rota and Sagan \cite{Rota}.

Recall that a \emph{disposition} is a function from $[m]$ to $[n]$ together with a
linear order on the pre-images of each $i\in [n]$.
Intuitively, a disposition can be visualized as a way of placing $m$
distinguished balls into $n$ distinguished boxes, where
the balls in each box are linearly ordered,
or equivalently, we may consider a disposition
as a decomposition of a permutation on $[m]$ into $n$ segments,
where we allow a segment to be empty.
We denote by $\mathcal{D}_{m,n}$ the set
of all dispositions from $[m]$ to $[n]$.

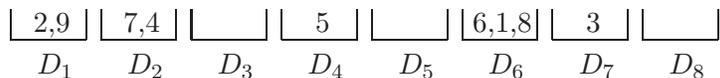
\begin{figure}[h,t]
\begin{center}
\begin{picture}(450,50)
\setlength{\unitlength}{1mm} \hspace{1cm}
\multiput(32,5)(12,0){8}{\line(1,0){10}}
\multiput(32,5)(12,0){8}{\line(0,1){4}}
\multiput(42,5)(12,0){8}{\line(0,1){4}}

\put(35,6){\small{2,9}}\put(47,6){\small{7,4}}
\put(72.5,6){\small{5}}\put(93.4,6){\small{6,1,8}}
\put(108.5,6){\small{3}}

\put(35.5,0.5){\small{$D_1$}} \put(47.5,0.5){\small{$D_2$}}
\put(59.5,0.5){\small{$D_3$}}\put(71.5,0.5){\small{$D_4$}}
\put(83.5,0.5){\small{$D_5$}}\put(95.5,0.5){\small{$D_6$}}
\put(107.5,0.5){\small{$D_7$}}\put(119.5,0.5){\small{$D_8$}}

\end{picture}
\vspace{0.1cm}\caption{An example of a disposition.}\label{fig7}
\end{center}
\end{figure}

For example,  Figure \ref{fig7} gives a disposition
from $[9]$ to $[8]$.

Let $D$ be a disposition from $[m]$ to $[n]$. We may write $D$ as
$(D_1,D_2,\ldots,D_n)$, where $D_1D_2\cdots D_n$ is a permutation of
$[m]$. Recall that, for a  permutation $\pi=\pi_1\pi_2\cdots \pi_k$
of $k$ elements, $\pi_i$ is said to be a \emph{right-to-left minimum} if
$\pi_i<\pi_j$ for each $j>i$. We denote by ${\rm RLmin}(D_i)$ the number
of right-to-left minima in $D_i$. For the disposition in Figure
\ref{fig7},  we have ${\rm RLmin}(D_1)=2$, ${\rm RLmin}(D_2)=1$,
${\rm RLmin}(D_3)=0$, ${\rm RLmin}(D_4)=1$, ${\rm RLmin}(D_5)=0$,
${\rm RLmin}(D_6)=2$, ${\rm RLmin}(D_7)=1$, ${\rm RLmin}(D_8)=0$.

As will be seen, the disposition polynomials are the generating
functions of   dispositions with respect to the statistics ${\rm
RLmin}(D_i)$. The proof of the following theorem is essentially
the same argument for the combinatorial interpretation of
the rising factorials.

\begin{theo}\label{thm3}
For $n\geq 1$, we have
\begin{equation}{\label{d1}}
\sum_{D\in \mathcal{D}_{m,n}}\prod_{i=1}^n
x_i^{\emph{RLmin}(D_i)}=\prod_{k=0}^{m-1}(x_1+x_2+\cdots+x_n+k).
\end{equation}
\end{theo}

\pf We use induction on $m$. For $m=1$, the assertion is
clear. Assume that (\ref{d1}) holds for $m-1$, that is,
\begin{equation}{\label{d2}}
\sum_{D\in \mathcal{D}_{m-1,n}}\prod_{i=1}^n
x_i^{{\rm RLmin}(D_i)}=\prod_{k=0}^{m-2}(x_1+x_2+\cdots+x_n+k).
\end{equation}

We proceed to show that the theorem holds for $m$.
A disposition from $[m]$ to $[n]$
can be obtained by inserting the element
 $m$ in a segment of a disposition from $[m-1]$ to $[n]$.
Let $(D_1, D_2,\ldots, D_n)$ be a disposition
from $[m-1]$ to $[n]$. Write $D_i=a_{1}a_{2}\cdots a_{r_i}$.
There are $r_i+1$ possible positions for the insertion of $m$ in $D_i$.
We consider two cases. Case 1. The element $m$ is attached to the
end of $D_i$. Let $D_i'=a_{1}a_{2}\cdots a_{r_i}m$. It is clear that $D_i'$ has one more right-to-left minima than $D_i$, that is,
\[{\rm RLmin}(D'_i)={\rm{RLmin}}(D_i)+1.\]
Case 2. The element $m$ is inserted before an element in $D_i$. Let $D_i'=
 a_{1}a_{2}\cdots a_{t-1}ma_{t}\cdots a_{r_i}$, for $1\leq t\leq r_i$. In this
case, we have
\[{\rm RLmin}(D'_i)={\rm{RLmin}}(D_i).\]
Since $r_1+r_2+\cdots +r_n=m-1$, considering
all possible insertions of $m$ into $(D_1,D_2,\ldots,D_n)$,
 we obtain that
\begin{align*}{\label{d3}}
\sum_{D\in \mathcal{D}_{m,n}}\prod_{i=1}^n
x_i^{{\rm{RLmin}}(D_i)}&=(x_1+r_1+\cdots +x_n+r_n) \sum_{D\in \mathcal{D}_{m-1,n}}\prod_{i=1}^n
x_i^{{\rm{RLmin}}(D_i)}\\
&=(x_1+x_2+\cdots+x_n+m-1)\sum_{D\in \mathcal{D}_{m-1,n}}\prod_{i=1}^n
x_i^{{\rm{RLmin}}(D_i)}.
\end{align*}
Thus, by the induction hypothesis,
we find that the theorem holds for $m$. This
completes the proof.\qed

In fact, one can  use the combinatorial
interpretation of the coefficients of the
rising factorials and the fundamental bijection for
permutations to deduce the above explanation of the
disposition polynomials. Recall that
the coefficient of $x^k$ in $(x)^m=x(x+1)\cdots(x+m-1)$
is the number of permutations of $[m]$ with $k$ cycles,
see Stanley \cite[1.3.4
Proposition]{Stanley1}.
The fundamental
bijection is also called the standard representation of
a permutation, or the first fundamental transformation.
For the purpose of this paper,
the standard representation of a permutation is defined as follows.
Based on the cycle decomposition, we write each cycle
by putting the minimum element at the end, and we arrange
the cycles  in the
increasing order of their minimum elements. Then we erase
all the parentheses.

Consider the set $\mathcal{S}_{m,n}$ of
cycle representations of permutations on $[m]$ with
each cycle colored by one of the $n$ colors, say $1, 2, \ldots, n$. For $\pi\in \mathcal{S}_{m,n}$, let ${\rm c}_{i}(\pi)$ denote the number of cycles of $\pi$ with
color $i$.

\begin{theo}\label{thmA1}
For $n\geq 1$, we have
\begin{equation}{\label{da1}}
\sum_{\pi\in \mathcal{S}_{m,n}}\prod_{i=1}^n
x_i^{{\rm c}_{i}(\pi)}=\prod_{k=0}^{m-1}(x_1+x_2+\cdots+x_n+k).
\end{equation}
\end{theo}

It can be seen that Theorem \ref{thm3} can be deduced from Theorem
\ref{thmA1} through the correspondence between permutations
with colored cycles and dispositions. For any $\pi\in
\mathcal{S}_{m,n}$, one may construct a
disposition $(D_1, D_2, \ldots, D_n)$  by the fundamental bijection, where $D_i$ is obtained
from the cycles of $\pi$ with color $i$.
Clearly, we have
\[{\rm RLmin}(D_i)={\rm c}_{i}(\pi).\]
Hence, Theorem \ref{thm3} can be deduced from Theorem \ref{thmA1}.

We define the homogenous disposition polynomials as follows
\[ Q_m(x_1,x_2,\ldots,x_n,t) = \prod_{k=1}^{m-1}(x_1+x_2+\cdots+x_n+kt).\]
For $n=m-1$, the homogenous disposition polynomials have been
 used by  Guo and Zeng \cite{GZ}.
Given a permutation $\pi=\pi_1\pi_2\cdots\pi_m$, Guo and Zeng defined
a general descent as an index  $i$  such that
$\pi_i>\pi_j$ for some $j>i$. Let
${\rm gdes}(\pi)$ denote  the number of general descents of $\pi$.
For a disposition
$D=(D_1,D_2,\ldots,D_n)$ from $[m]$ to $[n]$,
 let ${\rm gdes}(D)$ denote the
 total number of general descents of $D_i$ for $1\leq i \leq n$.
  It is easily checked  that
  \[{\rm
gdes}(D)=m-\sum_{i=1}^n{\rm RLmin}(D_i).\]
Then the homogeneous disposition polynomials have the following
combinatorial interpretation
 \begin{equation}
Q_m(x_1,x_2,\ldots,x_n,t)=\sum_{D\in \mathcal{D}_{m,n}}t^{{\rm
gdes}(D)}\prod_{i=1}^n x_i^{{\rm RLmin}(D_i)}.
\end{equation}

\section{A bijection between plane trees and dispositions}

In this section, we present a bijection
between plane trees and dispositions in the spirit of the
Pr\"ufer correspondence, which leads to a combinatorial
interpretation of the following formula of Guo and Zeng,
\begin{equation*}\label{e31}
\sum_{T\in \mathcal{P}_n}t^{{\rm eld}(T)}\prod_{i=1}^{n}x_i^{{\rm
young}_T(i)}=\prod_{k=0}^{n-2}(x_1+x_2+\cdots+x_n+kt),
\end{equation*}
where
$\mathcal{P}_{n}$ denotes the set of plane trees on $[n]$,
${\rm eld}(T)$ denotes the number of elder vertices of $T$ and ${\rm young}_T(i)$ denotes the number of younger children of vertex $i$ of $T$.

We now recall some terminology.
Given two vertices $i$ and $j$ of a plane tree $T$, we say that $j$
is a \emph{descendant} of $i$ if $i$ lies on the unique path from
the root to $j$. In particular, each vertex is a descendant of
itself. Denote by $\beta_T(i)$ the smallest descendant of $i$.
A child of $i$ means a descent $j$ such that $(i,j)$ is an edge of $T$.
A vertex $i$ is called the father of a vertex $j$ if $j$ is a child of
$i$. The vertices with the same father are called brothers of each other.
A vertex $j$ of a plane tree $T$ is called an elder vertex if $j$ has a
brother $k$ to its right such that $\beta_T(k)<\beta_T(j)$,
otherwise $j$ is called a younger vertex. Denote by ${\rm eld}_{T}(v)$
the number of elder children of $v$ in $T$, and denote by
${\rm young}_{T}(v)$ the number of younger children of $v$ in $T$.
It is not difficult to see that
${\rm young}_{T}(v)$ equals the number of right-to-left minima of the
sequence $\{\beta_T(v_1),\beta_{T}(v_2),\ldots,\beta_{T}(v_m)\}$,
where ${v_1,v_2,\ldots,v_m}$ are  the children of $v$ in linear order.
Moreover,
we denote by  ${\rm eld}(T)$ the total
number of elder vertices of $T$ and denote by ${\rm young}(T)$ the total number of younger  vertices of $T$.

\begin{figure}[h,t]
\begin{center}
\begin{picture}(130,50)
\setlength{\unitlength}{1mm}

\put(24,14){\hfill \rule{6pt}{6pt}}\put(25,18){\small$8$}
\put(25,5){\circle*{2}}\put(27,5){\small5}
\put(25,5){\line(0,1){10}}
\put(15,5){\circle*{2}}\put(12,5){\small2}
\put(15,5){\line(1,1){10}}
\put(34,4){\hfill \rule{6pt}{6pt}}\put(37,5){\small3}
\put(35,5){\line(-1,1){10}}
\put(5,-5){\circle*{2}}\put(0,-5){\small14}
\put(5,-5){\line(1,1){10}}
\put(24,-6){\hfill \rule{6pt}{6pt}}\put(27,-5){\small12}
\put(25,-5){\line(-1,1){10}}
\put(4,-16){\hfill \rule{6pt}{6pt}}\put(3.3,-19){\small16}
\put(5,-15){\line(0,1){10}}
\put(34,-6){\hfill \rule{6pt}{6pt}}\put(37,-4){\small1}
\put(35,-5){\line(0,1){10}}
\put(14,-16){\hfill \rule{6pt}{6pt}}\put(14.5,-19){\small4}
\put(15,-15){\line(2,1){20}}
\put(22,-16){\hfill \rule{6pt}{6pt}}\put(25,-18){\small6}
\put(23,-15){\line(6,5){12}}
\put(31,-15){\circle*{2}}\put(29.5,-19){\small 11}
\put(31,-15){\line(2,5){4}}
\put(38,-16){\hfill \rule{6pt}{6pt}}\put(38.5,-19){\small9}
\put(39,-15){\line(-2,5){4}}
\put(47,-15){\circle*{2}}\put(45.3,-19){\small 17}
\put(47,-15){\line(-6,5){12}}
\put(54,-16){\hfill \rule{6pt}{6pt}}\put(53.5,-19){\small13}
\put(55,-15){\line(-2,1){20}}
\put(22,-26){\hfill \rule{6pt}{6pt}}\put(21.3,-29){\small15}
\put(23,-25){\line(0,1){10}}
\put(33,-25){\circle*{2}}\put(31,-29){\small10}
\put(33,-25){\line(3,5){6}}
\put(44,-26){\hfill \rule{6pt}{6pt}}\put(44.5,-29){\small7}
\put(45,-25){\line(-3,5){6}}

\end{picture}
\vspace{3cm} \caption{A plane tree $T\in
\mathcal{P}_{17}$.}\label{fig5}
\end{center}
\end{figure}
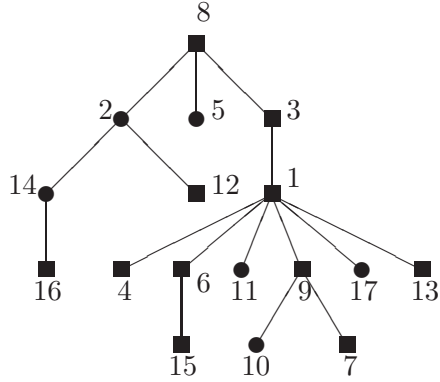

For
example, in Figure \ref{fig5}, each younger vertex is represented by
a  square, whereas each
elder vertex is represented by a solid dot.

\begin{theo}\label{thm0} There is a bijection $\varphi$
between plane trees on $[n]$ and
dispositions from $[n-1]$ to $[n]$. Let $T$ be a plane tree in $\mathcal{P}_n$,
and let $D=(D_1, D_2,\ldots,D_n)$ be the corresponding disposition under the
bijection $\varphi$. Then we have ${\rm young}_{T}(i)={\rm RLmin}(D_i)$ for all $i$.
\end{theo}
\pf We first give a description of the map $\varphi$ from $\mathcal{P}_n$ to $\mathcal{D}_{n-1,n}$. Let $T$ be a plane tree in $\mathcal{P}_n$.
We proceed to construct a disposition $D=(D_1,D_2,\ldots,D_n)$ through the
following procedure.

\begin{figure}[h,t]
\begin{center}
\begin{picture}(130,50)
\setlength{\unitlength}{1.2mm}

\put(25,15){\circle*{2}}\put(25,18){\small$8_0$}
\put(25,5){\circle*{2}}\put(27,5){\small$5_5$}
\put(25,5){\line(0,1){10}}
\put(15,5){\circle*{2}}\put(11,5){\small$2_3$}
\put(15,5){\line(1,1){10}}
\put(35,5){\circle*{2}}\put(37,5){\small$3_1$}
\put(35,5){\line(-1,1){10}}
\put(5,-5){\circle*{2}}\put(-2,-5){\small$14_{13}$}
\put(5,-5){\line(1,1){10}}
\put(25,-5){\circle*{2}}\put(27,-5){\small$12_{11}$}
\put(25,-5){\line(-1,1){10}}
\put(5,-15){\circle*{2}}\put(3.3,-19){\small$16_{15}$}
\put(5,-15){\line(0,1){10}}
\put(35,-5){\circle*{2}}\put(37,-4){\small$1_2$}
\put(35,-5){\line(0,1){10}}
\put(15,-15){\circle*{2}}\put(14.5,-19){\small$4_4$}
\put(15,-15){\line(2,1){20}}
\put(23,-15){\circle*{2}}\put(24,-18){\small$6_6$}
\put(23,-15){\line(6,5){12}}
\put(31,-15){\circle*{2}}\put(29.5,-19){\small $11_{10}$}
\put(31,-15){\line(2,5){4}}
\put(39,-15){\circle*{2}}\put(38.5,-19){\small$9_7$}
\put(39,-15){\line(-2,5){4}}
\put(47,-15){\circle*{2}}\put(45.3,-19){\small $17_{16}$}
\put(47,-15){\line(-6,5){12}}
\put(55,-15){\circle*{2}}\put(53.5,-19){\small $13_{12}$}
\put(55,-15){\line(-2,1){20}}
\put(23,-25){\circle*{2}}\put(21.3,-29){\small $15_{14}$}
\put(23,-25){\line(0,1){10}}
\put(33,-25){\circle*{2}}\put(31,-29){\small $10_9$}
\put(33,-25){\line(3,5){6}}
\put(45,-25){\circle*{2}}\put(44.5,-29){\small $7_8$}
\put(45,-25){\line(-3,5){6}}

\end{picture}
\vspace{3.5 cm} \caption{A plane tree with Pr\"ufer marks $T\in
\mathcal{P}_{17}$.}\label{fig4}
\end{center}
\end{figure}
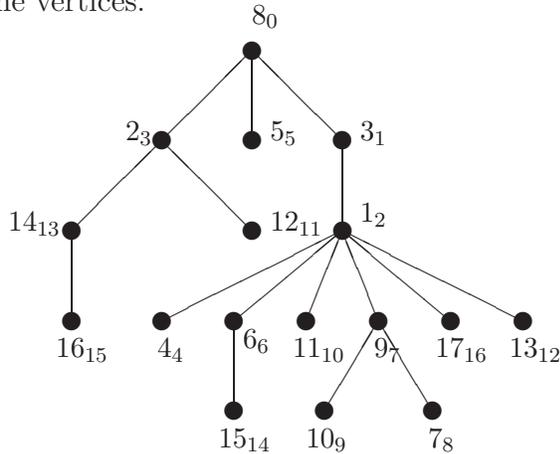

First, we mark the vertices of $T$ according to the
Pr\"ufer correspondence. More precisely, we mark the
vertices of $T$ by  the numbers
$0, 1, 2,\ldots, n-1$.  As the first step, we find the maximum leaf
of $T$, and mark it by $n-1$. Then we remove the maximum leaf and
repeat the this procedure until the root is marked by 0. These marks
are called the Pr\"ufer marks of $T$, which
represent the order that the vertices are removed in the Pr\"ufer
correspondence.
For example, Figure \ref{fig4} gives the Pr\"ufer marks
of a plane tree expressed by the indices of the vertices.

Using the Pr\"ufer marks, the disposition $D=(D_1, D_2,\ldots, D_n)$
can be easily constructed by setting $D_i$ to be the set of
the Pr\"ufer marks of the children of vertex $i$ endowed with the
linear order as in $T$. For example, for the plane tree $T$ in
Figure \ref{fig4},  we have $D_1=\{4,6,10,7,16,12\}$, $D_2=\{13,11\}$, $D_3=\{2\}$, and so on.

The above map $\varphi$ is indeed a bijection.
The inverse map can be described as follows. To
recover a plane tree $T$ from a disposition $D$,
we first mark the elements of $[n]$ by $0,1,2,\ldots,n-1$ from which
one recovers the Pr\"ufer marks of the plane tree $T$.
 We begin with the rightmost empty segment $D_i$, and mark the
 element $i$ by  $n-1$. Then we remove the empty segment $D_i$ and the element $n-1$ from some segment of $D$. Repeating this procedure until the last element of $[n]$ is marked by $0$.

\begin{figure}[h,t]
\begin{center}
\begin{picture}(450,50)
\setlength{\unitlength}{1mm}
\put(60,-40){\begin{picture}(450,50)
\put(16,15){\circle*{2}}\put(15.5,17){\small $2$}
\put(10,5){\circle*{2}}\put(7,4.5){\small $4$}
\put(10,5){\line(3,5){6}}

\put(22,5){\circle*{2}}\put(23.5,4.5){\small $5$}
\put(22,5){\line(-3,5){6}}

\put(10,-5){\circle*{2}}\put(9.5,-9){\small $6$}
\put(10,-5){\line(0,1){10}}

\put(16,-5){\circle*{2}}\put(15.5,-9){\small $3$}
\put(16,-5){\line(3,5){6}}

\put(28,-5){\circle*{2}}\put(27.5,-9){\small $1$}
\put(28,-5){\line(-3,5){6}}

\end{picture}
}

\put(-15,0){\begin{picture}(450,50)

\multiput(52,3)(14,0){6}{\line(1,0){10}}
\multiput(52,3)(14,0){6}{\line(0,1){4}}
\multiput(62,3)(14,0){6}{\line(0,1){4}}

\put(55,-1){\small{$D_1$}}\put(69,-1){\small{$D_2$}}
\put(83,-1){\small{$D_3$}} \put(97,-1){\small{$D_4$}}
\put(111,-1){\small{$D_5$}}\put(125,-1){\small{$D_6$}}

\put(69,4){\small{4,1}}\put(98,4){\small{5}}\put(111,4){\small{3,2}}
\end{picture}
}

\put(75,-12){$\Downarrow$}\put(80,-12){$\{6_5,4_4,3_3,1_2,5_1,2_0\}$}
\end{picture}
\vspace{4.5cm}\caption{An example to illustrate $\varphi^{-1}$ of the
case $n=6$.}\label{fig2}
\end{center}
\end{figure}
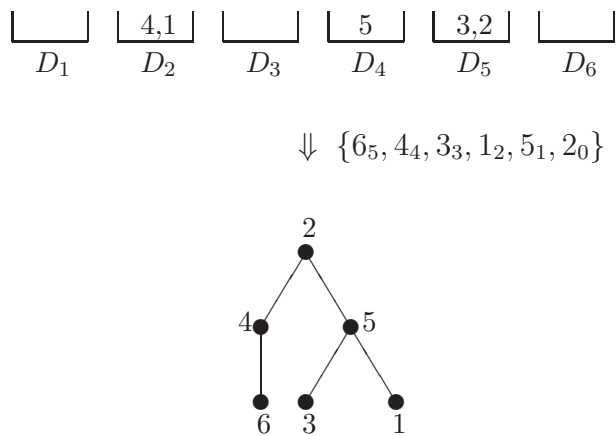

For example, for the disposition in Figure \ref{fig2}, the rightmost empty segment is $D_6$, thus, we mark $6$ by $5$. Deleting $D_6$ and and removing
$5$ from $D_4$, we see that $D_4$ becomes the rightmost empty segment.
 So we mark $4$ by $4$. Repeating this procedure, we
 obtain the marks $\{6_5,4_4,3_3,1_2,5_1,2_0\}$, where the index of each element stands for its mark.

Using the marks, we may construct the plane tree $T$ by setting the
root to be the element $r$ marked  by 0. If $D_r$ is empty, then $r$ must
be $1$ and $T$
consists of the single vertex $1$. Otherwise, we assume that $D_r=a_1,a_2,\ldots,a_t$, and assume that $b_i$ is
marked by $a_i$.  Set the children of  $r$ in linear order
to be $b_1,b_2,\ldots,b_t$.  Repeat the above process with respect to
each element $b_i$  until we
arrive at a plane tree $T$ on $[n]$.

Take the construction of the tree in Figure \ref{fig2} as an example.
We have already known the marks correspondence $\{6_5,4_4,3_3,1_2,5_1,2_0\}$. Notice that the element 2 is marked by 0, which indicates that 2 is the root of the corresponding tree. The elements in $D_2$ are $4,1$, which are the marks of $4,5$. Thus, the children of $2$ are $4,5$ in linear order. Similarly, the element in $D_4$ is $5$, which is the mark of the element $6$, thus, the only child of $4$ is $6$. Continue this procedure, and we will get the corresponding tree as demonstrated by Figure \ref{fig2}.

Now we aim to show that the above map is indeed the inverse of $\varphi$. It suffices to prove that the marks obtained from the disposition $D$
 are the same as the Pr\"ufer marks obtained from the plane tree $T$. Observe that the largest leaf $l$ in a plane tree on $[n]$ is marked by $n-1$.
 On the other hand, $D_l$ must be the rightmost segment in the corresponding disposition, and so $l$ is  marked by $n-1$ as well.
 We may repeat this argument for the element marked by $n-2$, if there is
 any segment left in the disposition. Hence we reach the conclusion that
 we get the same marks from the disposition $D$ and from the plane tree $T$.

 Next we  verify the relation
\[{\rm young}_{T}(i)={\rm RLmin}(D_i),\]
where $D$ is the corresponding disposition under $\varphi$.
It is not difficult to see that the degree of vertex $i$ of $T$ equals
the number of elements of $D_i$ in the disposition $\varphi(T)$.
Moreover,
let $D_i=b_1 b_2\cdots b_m$ and let $v_1,
v_2,\ldots,v_m$ be  the children of $i$ of $T$ in linear order.
We claim that for $1\leq j<k\leq m$, $b_j<b_k$ if and only if
$\beta(v_j) < \beta(v_k)$.
This property follows from the
fact that the Pr\"ufer mark of a vertex is the smallest among all its
descendants. Hence we deduce
that the number of younger children of vertex $i$ of $T$ equals the
number of right-to-left minima of $D_i$ in $\varphi(T)$. This
completes the proof.  \qed

It is clear that Theorem \ref{thm0} gives a combinatorial
interpretation of the following relation
\begin{equation}\label{eq1}
\sum_{T\in \mathcal{P}_n}t^{{\rm eld}(T)}\prod_{i=1}^n
x_i^{{\rm young}_T(i)}=\sum_{D\in \mathcal{D}_{n-1,n}}t^{{\rm
gdes}(D)}\prod_{i=1}^n x_i^{{\rm RLmin}(D_i)}.
\end{equation}

Combining (\ref{eq1}) and the combinatorial interpretation of the
disposition polynomials, we obtain a combinatorial proof of  the relation (\ref{e10}), that is,
\[\sum_{T\in \mathcal{P}_n}t^{{\rm eld}(T)}\prod_{i=1}^{n}x_i^{{\rm
young}_T(i)}=\prod_{k=0}^{n-2}(x_1+x_2+\cdots+x_n+kt).\]

Moreover, it can be seen that our correspondence can be
restricted to plane trees with a specific root $r$.
More precisely, a disposition $D$ corresponds to a plane tree $T$ with root $r$
if and only if the element $1$ is contained in $D_r$.
This leads to a combinatorial interpretation of relation
(\ref{e11}).

To conclude this section, we remark that the correspondence
$\varphi$ is also valid for labeled rooted trees. In this case,
we disregard the linear order in each segment of a disposition.
 In other words, $\varphi$ becomes a correspondence
 between labeled rooted trees and decompositions of $[n-1]$ into $n$ components.
Under this correspondence, the empty sets in a decomposition
correspond to leaves of a labeled rooted tree, and more generally,
the cardinalities of $D_i$ correspond to the degrees of the rooted trees.

\section{The Gessel-Seo polynomials}

In this section, we use the correspondence between plane trees
and dispositions to give a combinatorial interpretation of
the following expansion of  the Gessel-Seo polynomial,
\begin{equation}\label{e41}
\sum_{T \in
\mathcal{P}_{n+1}^{(1)}}x^{{\rm{young}}_{T}(1)}(t-z)^{{\rm{eld}}(T)}
z^{n-{\rm{young}}_{T}(1)-{\rm{eld}}(T)}=x\prod_{k=1}^{n-1}(x+(n-k)z+kt),
\end{equation}
where
$\mathcal{P}_{n+1}^{(1)}$ denotes the set of plane trees on $[n+1]$ with root $1$.
 Guo and Zeng \cite{GZ} derived the above identity by using generating
 functions and asked for a combinatorial proof.
In fact, as a consequence of (\ref{e11}), that is,
\[\sum_{T\in \mathcal{P}_n^{(r)}}t^{{\rm eld}(T)}\prod_{i=1}^{n}x_i^{{\rm
young}_T(i)}=x_r\prod_{k=1}^{n-2}(x_1+x_2+\cdots+x_n+kt),\]
we find that (\ref{e41}) holds for plane trees on $[n+1]$ with any specific root $r$ by
replacing $n$ by $n+1$ and setting $x_r=x$, $x_i=z$ for any
for $i\not=r$. As will be seen, our correspondence between plane trees and
dispositions  serves as a direct combinatorial interpretation of
this fact.

\begin{theo}For $n\geq 1$ and $1\leq r \leq n+1$, we have
\begin{equation}\label{e4}
\sum_{T\in
\mathcal{P}_{n+1}^{(r)}}x^{\emph{young}_{T}(r)}(t-z)^{\emph{eld}(T)}
z^{n-\emph{young}_{T}(r)-\emph{eld}(T)}=x\prod_{k=1}^{n-1}(x+(n-k)z+kt),
\end{equation}
where $\mathcal{P}_{n+1}^{(r)}$
  is the set of plane trees on $[n+1]$ with root $r$.
\end{theo}

\pf Replacing $t$ by $t+z$, we may rewrite $(\ref{e4})$ as follows,
\begin{equation}\label{e5}
\sum_{T\in \mathcal{P}_{n+1}^{(r)}}x^{{\rm
young}_{T}(r)}t^{{\rm eld}(T)}
z^{n-{\rm young}_{T}(r)-{\rm eld}(T)}=x\prod_{k=1}^{n-1}(x+nz+kt).
\end{equation}

We first give a combinatorial interpretation of the right hand side of ($\ref{e5}$).
By the combinatorial interpretation of the disposition polynomials,
 we see that the Gessel-Seo polynomial $P_n(t+z,z,x)$ is the generating function
 of  dispositions $D=(D_1,D_2,\ldots,D_{n+1})$ from
$[n]$ to $[n+1]$ with the element $1$
contained in $D_r$, where a right-to-left minimum in $D_r$ is given a weight $x$, a right-to-left minimum in $D_i$ ($i\not=r$) is given a weight $z$, and any other element is
given a weight $t$.

For a disposition $D=(D_1,D_2,\ldots,D_{n+1})$ with which the element $1$ appears in $D_r$, let $T$ be the plane tree corresponding to $D$ under the bijection
 $\varphi$ in Theorem \ref{thm0}. It is easily seen that $T$ has root $r$, namely, $T\in \mathcal{P}_{n+1}^{(r)}$. Moreover, a younger child of vertex $i$ of $T$ corresponds to a right-to-left minimum in $D_i$, and an elder child of
vertex $i$ of  $T$ corresponds to an element which is not a right-to-left minimum in $D_i$ for $1\leq i\leq n+1$. Hence $T$ has weight
\[ x^{{\rm
young}_{T}(r)}t^{{\rm eld}(T)}
z^{n-{\rm young}_{T}(r)-{\rm eld}(T)},\]
as expected. This completes the proof.\qed

\vskip 5mm \noindent{\bf Acknowledgments.} This work was supported
by the 973 Project, the PCSIRT Project of the Ministry of Education,
and the National Science Foundation of China.


\begin{thebibliography}{99}

\bibitem{Gessel} I.M. Gessel and S. Seo, A refinement of Cayley's
formula for trees, Electron. J. Combin, 11 (2004), \# R27.

\bibitem{GZ} V.J.W. Guo and J. Zeng, A generalization of the
Ramanujan polynomials and plane trees, Adv. Appl. Math. 39 (2007),
96--115.

\bibitem{Rota} S.A. Joni, G.C. Rota and B. Sagan, From sets to
functions: Three elementary examples, Discrete Math. 37 (1981),
193--202.

\bibitem{Rota1} R. Mullin and G.C. Rota, On the foundations of
combinatorial theory III: The theory of binomial enumeration, Graph
Theory and Applications, Academic Press, New York, 1970, pp. 167--213.



\bibitem{Shin} S. Seo and H. Shin, A generalized enumeration of
labeled trees and reverse Pr\"ufer algorithm, J. Combin. Theory,
Ser. A, 114 (2007), 1357--1361.

\bibitem{Shin1} H. Shin, A new bijection between forests and parking
functions, arXiv: 0810.0427.

\bibitem{Stanley1} R.P. Stanley, Enumerative Combinatorics, Vol. 1,
Cambridge University Press, Cambridge, 1997.


\end{thebibliography}
\end{document}